\documentclass[preprint,review,12pt,authoryear]{elsarticle}

\usepackage{amssymb}

\usepackage{mathrsfs}
\usepackage[colorlinks, linkcolor = red]{hyperref}
\usepackage{xcolor}
\usepackage{mathscinet}
\usepackage{latexsym}
\usepackage{amsthm}
\usepackage{amssymb}
\usepackage{amsfonts}
\usepackage{amsmath}
\usepackage{longtable}
\usepackage{graphicx}
\usepackage{multirow}
\usepackage{multicol}
\usepackage{soul}
\usepackage{booktabs}
\usepackage{cases}
\usepackage{natbib}
\usepackage{enumerate}
\usepackage{fullpage}

\newtheorem{theorem}{Theorem}[section]
\newtheorem{thm}[theorem]{Theorem}

\newtheorem{lem}[theorem]{Lemma}
\newtheorem{proposition}[theorem]{Proposition}

\newtheorem{corollary}[theorem]{Corollary}

\theoremstyle{definition}

\newtheorem{defn}[theorem]{Definition}

\theoremstyle{remark}

\newtheorem{rem}[theorem]{Remark}
\numberwithin{equation}{section}

\DeclareMathAlphabet{\mathpzc}{OT1}{pzc}{m}{it}

\newcommand{\Be}{\begin{equation}}
\newcommand{\Ee}{\end{equation}}
\newcommand{\Bs}{\begin{split}}
	\newcommand{\Es}{\end{split}}
\newcommand{\Bes}{\begin{equation*}}
	\newcommand{\Ees}{\end{equation*}}
\newcommand{\BT}{\begin{thm}}
	\newcommand{\ET}{\end{thm}}
\newcommand{\Bp}{\begin{proof}}
	\newcommand{\Ep}{\end{proof}}
\newcommand{\BL}{\begin{lem}}
	\newcommand{\EL}{\end{lem}}
\newcommand{\BP}{\begin{proposition}}
	\newcommand{\EP}{\end{proposition}}
\newcommand{\BC}{\begin{corollary}}
	\newcommand{\EC}{\end{corollary}}
\newcommand{\BR}{\begin{rem}}
	\newcommand{\ER}{\end{rem}}
\newcommand{\BD}{\begin{defn}}
	\newcommand{\ED}{\end{defn}}
\newcommand{\BI}{\begin{itemize}}
	\newcommand{\EI}{\end{itemize}}

\newcommand{\design}{\boldsymbol{X}}
\newcommand{\designrow}{\boldsymbol{x}}
\newcommand{\designs}{x}

\newcommand{\outcomes}{y}

\newcommand{\parameter}{\boldsymbol{\beta}}
\newcommand{\target}{\parameter^*}
\newcommand{\est}{\widehat{\parameter}}
\newcommand{\parameters}{\beta}
\newcommand{\targets}{\parameters^*}

\newcommand{\difvec}{\boldsymbol{\delta}}
\newcommand{\nzeroindx}{\boldsymbol{S}}

\newcommand{\vecv}{\boldsymbol{v}}
\newcommand{\vecu}{\boldsymbol{u}}

\newcommand{\contx}{R}

\newcommand{\tp}{^{\top}}

\journal{}

\begin{document}

\begin{frontmatter}



\title{Consistency of $\ell _{1}$ Penalized Negative Binomial Regressions}


\author{Fang Xie\fnref{label1}}
\address{School of Mathematics and Statistics, Wuhan University,\\
	Wuhan, Hubei 430072, P.R. China.\\}
\ead{fangxie219@foxmail.com}
\fntext[label1]{Corresponding author}

\author{Zhijie Xiao}
\address{Department of Economics, Boston College, Chestnut Hill, MA 02467}
\ead{zhijie.xiao@bc.edu}


\begin{abstract}
We prove the consistency of the $\ell_1$ penalized negative binomial regression (NBR). 
A real data application about German health care demand shows that the $\ell_1$ penalized NBR produces a more concise but more accurate model, comparing to the classical NBR. 
\end{abstract}



\begin{keyword}


negative binomial regression \sep $\ell _{1}$ consistency \sep penalized maximum
likelihood \sep high-dimensional regression
\end{keyword}

\end{frontmatter}



\section{Introduction}

Count data is an important type of statistical data in which the observation takes non-negative integer values. 
Count data naturally arises in many
areas such as health care demand \citep{RWM03}, consumer credit
behaviors \citep{GWH94}, vehicle crash \citep{WL13}, psychology %
\citep{GMS95} and so on. 
Poisson distribution is a counting measure
extensively used to model count data \citep{CT98}, and the
Poisson regression has been an important generalized linear model that is
widely used in applications \citep{CDL15,SSL08}. 
Moreover, penalized Poisson regressions have been extensively studied and used to model high dimensional count data \citep{AL15,IPR16,LC15}. 
However, a major limitation of the Poisson regression is its restrictive assumption
that the variance equals the mean. 
In practice, more and more
applications are found to have an overdispersion feature that sample
variance is much larger than sample mean \citep{CT98,Hilbe11}, which
violates the assumptions of Poisson regression. 
For this reason, a
more general and flexible regression model, the negative binomial
regression, has attracted a great deal of research attention and become
a popular model in analyzing count data. 

The NBR, as a generalization of the Poisson regression, loosens the
highly restrictive assumption that the variance is equal to the mean made by
the Poisson model. 
The negative binomial distribution has two parameters, the mean parameter $\mu $, and the over-dispersion parameter $r$. 
dispersion property. 
When $r\rightarrow \infty $, the negative binomial
distribution converges to a Poisson distribution with the parameter $\mu $ %
\citep{CT98,Hilbe11}. 

Nowadays, negative binomial distribution is becoming more and more important in modeling real data in health care science \citep{LWLM13,RWM03}, biology \citep{MDS15}, psychology \citep{Wal07}, medicine \citep{ACH14,AWFPS16}, ecology \citep{LM11}, finance \citep{CT96} and so on. 
As the dimension of data increases, variable selection is very important and
necessary to simplify the fitting models. 
Stimulated by the great success of many penalized regressions such as lasso \citep{Tib96}, the NBRs with a penalty have recently been proposed to analyze high dimensional data, for example, the data of the association between multiple
biomarkers and prolonged hospital length of stay \citep{WMZPWD16}. 
However,
there are a few literature about the penalized negative binomial
regression and hence less statistical theories. 
So, the first and main goal of this paper is to rigorously prove the consistency property of $\ell_1$ penalized NBR. 


In addition to the theoretical analysis, we also apply the $\ell _{1}$ penalized NBR for analyzing real data about German health care demand. 
The data, supplied by the German Socioeconomic Panel (GSOEP),  
consist of 27326 samples observed in seven year. 
There are two dependent variables and 23 variables. 
For the brevity of analysis, we only consider one dependent variable, the
number of doctor visit within the last quarter prior to the survey (DOCVIS). 
We compare the $\ell _{1}$ penalized NBR with the classical NBR. 
The results show the $\ell _{1}$ penalized NBR can produce a more simple and efficient model, which exhibits the most effective variables and gives a much smaller prediction error. 


The rest of the paper is organized as follows. 
In Section \ref{ssec:notation}, we give the notations throughout the paper. 
In Section \ref{ssec:model}, we introduce the model and the $\ell _{1}$ penalized NBR method. 
The theoretical results are shown in Section \ref{ssec:thmres}.
Sections \ref{sec:experiments} and \ref{sec:realdata} show the numerical
results of $\ell _{1}$ penalized NBR based on synthetic data and real data. 
We finally conclude in Section~\ref{sec:conclusion}. 

\subsection{Notations}\label{ssec:notation}
For any vector $\boldsymbol{v}=(v_{1},\cdots
,v_{d})\in \mathbb{R}^{d}$, $\Vert \boldsymbol{v}\Vert _{q}$ denotes its $l_{q}$-norm
with $0<q<\infty $. When $q=\infty $, $\Vert \boldsymbol{v}\Vert _{\infty
}=\max_{i\in \lbrack d]}|v_{i}|$ where the notation $[d]=\{1,2,\cdots
,d\}$. 
Let $\{a_{n}\}_{n\geq 1},\{b_{n}\}_{n\geq 1}$ be two real sequences. The
notation $b_{n}=o(a_{n})$ means that $\lim_{n\rightarrow \infty }b_{n}/a_{n}=0$. The notation $b_{n}=O(a_{n})$ means that $\lim_{n\rightarrow
	\infty }b_{n}/a_{n}=C$ with $C$ being some constant.


\section{Model and Main theorems}\label{sec:main} \label{sec:thm}
In this section, we first introduce the negative binomial model and the corresponding $\ell_1$ regularized method. 
Then, we give our theoretical results on the consistency of this method. 

\subsection{Model and penalized negative binomial regression}

\label{ssec:model} The negative binomial regression model assumes that
the response variable $Y$ has a negative binomial distribution, and the
logarithm of its expected value can be modeled by a linear combination of
unknown parameters.\ Thus, given observed data $(\designrow_{1},\outcomes_{1})$, $\cdots
,(\designrow_{n},\outcomes_{n})$, we assume that 
\begin{equation*}
\outcomes_{i}|\designrow_{i} \sim \mathrm{NB}(r,\mu _{i})\ \mathrm{with}\
\mu _{i}=e^{\designrow_{i}\tp \target}\ \mathrm{and}\ r>0,
\end{equation*}%
with $\outcomes_{i}\in \mathbb{R}$, $\designrow_{i}\in \mathbb{R}^{p}$ and $\mathrm{NB}(r,\mu
_{i})$ signifies a negative binomial distribution with parameters $(r,\mu
_{i})$, where the mean $\mu _{i}$ is modeled as $e^{\designrow_{i}\tp \target}$%
, and $\target=(\targets_{1},\cdots ,\targets_{p})\tp \in 
\mathbb{R}^{p}$ is an unknown parameter to be estimated.

By definition, the negative binomial probability density function has the
form\citep{WMF11}
\begin{equation*}
\mathbb{P}(\outcomes_{i}|\designrow_{i})=\left[ \frac{r}{r+\mu _{i}}\right] ^{r}%
\frac{\Gamma (r+\outcomes_{i})}{\Gamma (r)\outcomes_{i}!}\left[ \frac{\mu _{i}}{r+\mu _{i}}%
\right] ^{\outcomes_{i}}\text{, }\mathrm{with}\ \mu _{i}=e^{\designrow_{i}\tp\target}%
\text{,}\ 
\end{equation*}%
$\mathrm{for\ }i=1,2,\cdots ,n$, where $\Gamma (\cdot )$ is the gamma
function.

Denote $\designrow_{i}=(\designs_{i1},\cdots
,\designs_{ip})\tp $. Without loss of generality, we assume that 
\begin{equation*}
\frac{1}{n}\sum\limits_{i=1}^{n}x_{ij}=0\quad \mathrm{and}\quad \frac{1}{n}%
\sum\limits_{i=1}^{n}x_{ij}^{2}=1.
\end{equation*}%

The estimator $\est$ for negative binomial regression obtained
by $\ell _{1}$ penalized maximum log-likelihood method is defined by 
\begin{equation}
\est=\arg \min\limits_{\parameter\in \mathbb{R}^{p}}\left\{ L(\parameter)+\lambda \Vert \parameter \Vert _{1}\right\}, \label{def:betahat}
\end{equation}
where 
\begin{equation*}
L(\parameter)=-\frac{1}{n}\sum\limits_{i=1}^{n}\left( \outcomes_{i}(\designrow_{i}\tp \parameter-\ln
(r+e^{\designrow_{i}\tp \parameter}))-r\ln (r+e^{\designrow_{i}\tp \parameter})\right) .
\label{eq:fbeta}
\end{equation*}%

From Karush-Kuhn-Tucker conditions, we know that $\lambda \geq \Vert
\nabla L(\est)\Vert _{\infty }$. So, if the $\est$ is close to 
$\target$, the event $\{\lambda \geq c\Vert \nabla
L(\target)\Vert _{\infty }\}$ will hold with high probability with some
constant $c>1$. In fact, we indeed prove this event holds with probability
approaching to 1 as $n,p\rightarrow \infty $, see Lemmas \ref{lem:asym} and %
\ref{lem:exact} below.

\subsection{Theoretical results}

\label{ssec:thmres} We aim to prove the consistency of estimator $\est
$ and obtain the convergence rate of estimation error.

 Firstly, we denote $V=\Vert \nabla L(\target)\Vert _{\infty }
$ where 
\begin{equation*}
\nabla L(\target)=-\frac{1}{n}\sum\limits_{i=1}^{n}\designrow_{i}\times \sqrt{%
	\frac{re^{\designrow_{i}\tp \target}}{e^{\designrow_{i}\tp \target}+r}}\times 
\frac{(\outcomes_{i}-e^{\designrow_{i}\tp \target})}{\sqrt{\frac{e^{\designrow_{i}\tp \parameter
				^{\ast }}(r+e^{\designrow_{i}\tp \target})}{r}}}.  \label{eq:gradient}
\end{equation*}%
Notice that the last term on the right-hand side of equation above is the normalized form of $%
\outcomes_{i}$. This part will play an important role later in evaluating the
probability of event $\{\lambda \geq c\Vert \nabla L(\target)\Vert
_{\infty }\}$. For each $1\leq i\leq n$, denote 
\begin{equation*}
v_{i}=\sqrt{\frac{re^{\designrow_{i}\tp \target}}{e^{\designrow_{i}\tp \target}+r%
}}\ \ \ \ \ \ \ \mathrm{and}\ \ \ \ \ \ \ v_{n}=\max\limits_{1\leq i\leq
	n}v_{i}.  \label{eq:vn}
\end{equation*}%
For the negative binomial distribution, it is easy to check that its mean and
variance have the following relationship: 
\begin{equation*}
\mathrm{Var}(\outcomes_{i}|\designrow_{i})=\mu _{i}\left( 1+\frac{\mu _{i}}{r}\right) ,
\end{equation*}%
where $\mu _{i}$ is the mean and the parameter $r$ describes the dispersion
of the negative binomial distribution. If there exists some positive constant $%
B>1$ such that $r\leq \mu _{i}/B-1$ for all $i$, we obtain that $%
\mathrm{Var}(\outcomes_{i}|\designrow_{i})/\mu _{i}\in \lbrack B,+\infty )$ for
all $i$. Hence, $v_{i}\leq \sqrt{\mu _{i}/B}$ for each $i$ and $%
v_{n}\leq \max_{1\leq i\leq n}\sqrt{\mu _{i}/B}$. This
implies that the value of $v_{n}$ depends on the maximum of means $\mu _{i}$.%

For notation ease, we itemize the conditions throughout the paper:
\begin{itemize}
	\item[\textbf{C1}.] Sparsity of $\target$: $s=|\nzeroindx|<n$ where $\nzeroindx=\{j\in[p]:
	\targets_{j}\neq 0\}$.
	\item[\textbf{C2}.] There exists a positive constant $r$ such that $%
	\sup_{i\in[n],j\in[p]}|x_{ij}|\le \contx<\infty$.
	\item[\textbf{C3}.] $n,p$ satisfy that $\sqrt{n}< p\le o(e^{n^{1/5}})$ and $%
	p/\alpha>8$ for all $\alpha\in(0,1)$.
	\item[\textbf{C4}.] For any $\difvec\in\mathbb{R}^p$ satisfying $\|\difvec_{\nzeroindx^{c}}\|_1\le \gamma\|\difvec_{\nzeroindx} \|_1$ with some $
	\gamma>1$, there exists a positive constant $\phi_0$ such that $%
	\langle \difvec, \nabla^2 L(\target)\difvec \rangle\ge\phi_0^2\|\difvec_{\nzeroindx} \|_2^2$.
\end{itemize}

 Conditions \textbf{C1} and \textbf{C2} are conventional conditions
that have been widely used in literature. 
Condition \textbf{C3} is assumed since we focus on the high-dimensional regression problems including $n>p$
and $n<p$. 
Condition \textbf{C4} is the popular restricted eigenvalue
assumption, see also in \citep{BRT09,JXX17}.

 The following theorem gives the bounds of $|L(\est%
)-L(\target)|$ and $\ell_1$ estimation error of $\est$. 

\begin{thm}
	\label{thm:main} Let $\est$ be defined in \eqref{def:betahat}.
	Suppose that the conditions \textbf{C1}, \textbf{C2}, and \textbf{C4} hold.
	If $\lambda\ge cV$ with some $c>1$ and $\lambda s\le(c-1)^{2}\phi_0^{2}/6c\contx(c+1)$ with $s, \contx, \phi_{0}$ defined in \textbf{C1}, \textbf{C2}, and 
	\textbf{C4} respectively, we have 
	\begin{align}  
	\|\est-\target\|_1&\le\frac{C\lambda s}{\phi_0^2},\label{ineq:main1} \\
	|L(\est)-L(\target)|&\le\frac{C\lambda^2 s}{\phi_0^2}, \label{ineq:main2}
	\end{align}
	where $C=2c_{1}c(c+1)/(c-1)^{2}$ with some constant $c_{1}\in(2,3]$.
\end{thm}

The proof of Theorem \ref{thm:main} is given in the Supplementary Material (see Appendix A).
From inequalities \eqref{ineq:main1} and \eqref{ineq:main2}, we can know
that $\|\est-\target\|_1=O(\lambda s)$ and $|L(\est%
)-L(\target)|=O(\lambda^2 s)$ which are satisfactory. But these two orders
are obtained based on the event $\{\lambda\ge cV\}$ with some $c>1$ and the
condition $\lambda s\le(c-1)^{2}\phi_0^{2}/6c\contx(c+1)$. So, we need to
have further studies on the event and the order of $s$. The order of $s$ is
discussed in Remark \ref{rem:order}.

We first study the event $\{\lambda \geq cV\}$ with some $c>1$. Recall $%
V=\Vert \nabla L(\target)\Vert _{\infty }$ and denote by $V(1-\alpha |%
\design)$ the $(1-\alpha )$-quantile of $V$. We give two choices of $\lambda $
as follows. Given any $\alpha \in (0,1)$ and some $c>1$, define 
\begin{align*}
\mathrm{exact\ choice:}\ \lambda & =cV(1-\alpha |\design), \\
\mathrm{asymptotic\ choice:}\ \lambda & =cv_{n}(\sqrt{n})^{-1}\Phi ^{-1}(1-%
\frac{\alpha }{2p}),  
\end{align*}%
where $v_{n}$ is defined in \eqref{eq:vn} and $\Phi $ is the cumulative
distribution function of standard normal distribution.

We will show in the following two lemmas that under these two choices of $%
\lambda$, the probabilities of the event $\{\lambda\ge cV\}$ will approach
to 1 as $n,p\rightarrow\infty$.

\begin{lem}
	\label{lem:exact} If $\lambda=cV(1-\alpha|\design)$ with some $\alpha\in
	(0,1) $, then we have 
	\begin{equation*}
	\mathbb{P}(\lambda\ge cV)\ge 1-\alpha.
	\end{equation*}
\end{lem}
Lemma \ref{lem:exact} can be easily proved by the definition of the quantile.

\begin{lem}
	\label{lem:asym} If $\lambda=cv_n(\sqrt{n})^{-1}\Phi^{-1}(1-\alpha/2p%
	)$ with some $\alpha\in (0,1)$, then we have 
	\begin{align*}
	\mathbb{P}(\lambda\ge cV)&\ge 1-\alpha\left(1+O(1)(\sqrt{2\log{(2p/\alpha)}}-%
	\sqrt{n}m)^{3}n^{-1/2}(3w_{1}\log{p}+m)\right) \\
	&\quad\times\left(1+\frac1{\log{(p/\alpha)}}\right)\frac{\exp\{(n\log{(p/\alpha)}%
		)^{1/2}m-nm^{2}/2\}}{1-\sqrt{n}m/(\log{(p/\alpha)})^{1/2}}+C_{1}n/p^{2},
	\end{align*}
	where $m=6C_{1}w_{1}\log{p}/p^{3}$ with some positive constants $%
	C_{1}$ and $w_{1}$. In particular, as $n,p\rightarrow\infty$, we have 
	\begin{equation*}
	\mathbb{P}(\lambda\ge cV)\ge 1-\alpha(1+o(1)).
	\end{equation*}
\end{lem}

The proof of Lemma \ref{lem:asym} is given in the Supplementary Material (see Appendix B).
The key technique is Cram{\'e}r type moderate deviation theorem. 

Combing Lemmas \ref{lem:exact} and \ref{lem:asym} with
Theorem \ref{thm:main}, we can obtain the following propositions \ref{prop:exact} and \ref{prop:asym} respectively, which give the bounds of $\ell_{1}$ estimation error under two choices of $\lambda$. 

\begin{proposition}
	\label{prop:exact} Let $\est$ be defined in \eqref{def:betahat}.
	Suppose that the conditions \textbf{C1}, \textbf{C2}, \textbf{C3}, and 
	\textbf{C4} hold. If $\lambda=cV(1-\alpha|\design)$ with some $\alpha\in
	(0,1) $ and $c>1$ and $\lambda s\le(c-1)^{2}\phi_0^{2}6c\contx(c+1)$
	with $s, \contx, \phi_{0}$ defined in \textbf{C1}, \textbf{C2}, and \textbf{C4},
	then with probability at least $1-\alpha$, the inequalities~\eqref{ineq:main1} and~\eqref{ineq:main2} hold. 
\end{proposition}

\begin{proposition}
	\label{prop:asym} Let $\est$ be defined in \eqref{def:betahat}. Suppose
	that the conditions \textbf{C1}, \textbf{C2}, \textbf{C3}, and \textbf{C4}
	hold. If $\lambda=cv_n(\sqrt{n})^{-1}\Phi^{-1}(1-\alpha/2p)$ with
	some $\alpha\in (0,1)$ and $c>1$ and $\lambda s\le (c-1)^{2}\phi_0^{2}/6c\contx(c+1)$ with $s, \contx, \phi_{0}$ defined in \textbf{C1}, \textbf{C2}, and 
	\textbf{C4}, then with probability at least 
	\begin{align*}
	& 1-\alpha\left(1+O(1)(\sqrt{2\log{(2p/\alpha)}}-\sqrt{n}%
	m)^{3}n^{-1/2}(3w_{1}\log{p}+m)\right) \\
	&\quad\times(1+\frac1{\log{(p/\alpha)}})\frac{\exp\{-2(n\log{(p/\alpha)}%
		)^{1/2}m+nm^{2}\}}{1-\sqrt{n}m/(\log{(p/\alpha)})^{1/2}}+C_{1}n/p^{2},
	\end{align*}
	with $m=6C_{1}w_{1}\log{p}/p^{3}$ with some positive constants $%
	C_{1} $ and $w_{1}$, the inequalities~\eqref{ineq:main1} and~\eqref{ineq:main2} hold. 
\end{proposition}

\begin{rem}
	\label{rem:order} Notice that the asymptotic choice of $\lambda$ is order of 
	$\sqrt{(\log{p})/n}$. Then, $\|\est-\target\|_1=O(s\sqrt{(\log{p})/n})$ and $|L(\est)-L(\target)|=O(s(\log{p})/n)$. If $s$ satisfies $s=o(\sqrt{n/\log{p}})$, then $\|\est-\target\|_1=o(1)$ and $|L(\est)-L(\target)|=o(\sqrt{(\log{p})/n})$.
\end{rem}


\section{Simulations}

\label{sec:experiments} In this section, we conduct simulations to show the
performances of $\ell_{1}$ penalized NBR in the perspective of the variation of the estimation error of estimator $\est$, true positive rate and true negative rate.  %
%

We generated the data $(\outcomes_i,\designrow_i)$ from the negative binomial distribution $%
\mathrm{NB}(r,\mu_i)$ with the following setting: $r$ varies in $\{2,1,0.5,0.25\}$, $\mu_i=e^{\designrow_{i}\tp\target}$ for $i\in[n]$, $p$-dimensional observations $\designrow_{i}\sim \mathrm{N}(0,\boldsymbol{\Sigma})$ with $\Sigma_{jk}=\rho^{|j-k|}$ for $j,k\in[p]$,
$p$-dimensional true parameter vector $\target$
has $5$ nonzero components taking value in $[-1, 1]$ randomly.
We set $p=30$, $n\in\{100,200,400,800\}$, and $\rho=0.5$. The smaller $r$ means that the data $\outcomes_{i}$ are more over-disperse. 

We will compare the estimation error of estimator $\est$, true positive rate
(sensitivity) and true negative rate (specificity). 
We define the estimation error by $\|\est -\target\|_1$. Sensitivity is the fraction of the number of correctly selected predictors in all the effective predictors, and specificity is the fraction of the number of correctly unselected predictions in all the ineffective predictors.  
An ideal estimator should have the estimation error close to 0, and 
sensitivity and specificity close to 1. 
We repeat each realization 100 times to obtain the mean and standard deviation (SD) of the three criteria above. 
The results are shown in Table \ref{tab:sim1}. It indicates that for fixed $(r,p,\rho)$, the estimation error decreases with $n$ increasing. At
the same time, sensitivity increases towards 1 as we expect, but this pays a price on specificity. So, we could choose an
appropriate sample size to make a trade-off between sensitivity and
specificity based on the actual situation.
\begin{table}[tbp]
	\centering
	\scriptsize
	\caption{Estimation error, sensitivity, specificity of negative binomial regression with $\ell_{1}$ penalty under different cases $(n,p,r,\protect\rho)$.}
	\label{tab:sim1}%
\begin{tabular}{cccccccc}
	\hline
	\multirow{2}{*}{$n$}& \multicolumn{3}{c}{$r=2, p = 30, \rho = 0.5$} && \multicolumn{3}{c}{$r=1, p = 30, \rho = 0.5$}\\
	& estimation error (SD) & sensitivity & specificity && estimation error (SD) & sensitivity & specificity \\
	\hline
	100	&   0.484(0.188)    &	 0.716	&	0.702	&&	0.516(0.182)	&	0.488	&	0.798		\\
	200	&	0.245(0.123)    &	0.848	&	0.653	&&	0.332(0.121)	&	0.628	&	0.766		\\
	400	&	0.121(0.061)	&	0.924	&	0.605	&&	0.184(0.096)	&	0.836	&	0.670		\\
	800	&	0.061(0.033)	&	0.984	&	0.555	&&	0.093(0.032)	&	0.980	&	0.573		\\
	\hline
	\\
	& \multicolumn{3}{c}{$r=0.5, p = 30, \rho = 0.5$} && \multicolumn{3}{c}{$r=0.25, p = 30, \rho = 0.5$}\\
	& estimation error (SD) & sensitivity & specificity && estimation error (SD) & sensitivity & specificity \\
	\hline
	100	&	0.785(0.363)	&	0.488	&	0.787	&&	1.081(0.638)	&	0.284	&	0.826		\\
	200	&	0.480(0.210)	&	0.640	&	0.734	&&	0.634(0.273)	&	0.560	&	0.795		\\
	400	&	0.241(0.126)	&	0.824	&	0.677	&&	0.392(0.185)	&	0.704	&	0.759		\\
	800	&	0.122(0.032)	&	0.940	&	0.603	&&	0.210(0.121)	&	0.848	&	0.635		\\
	\hline
\end{tabular}
\end{table}

%
%

\section{An Application on German health care demand}
\label{sec:realdata} 
We will apply the NBR with the $\ell _{1}$-penalized MLE
method to a real dataset on German health care demand in this section. 
The data is a part of the German Socioeconomic Panel
(GSOEP) data which was employed in \citep{RWM03}. The data source can
be downloaded on {\url{http://qed.econ.queensu.ca/jae/2003-v18.4/riphahn-wambach-million/}}.
The data consist of 27326 observations from 7293 individuals observed one or
several times during years $\{1984, 1985,$ $1986,1987, 1988, 1991,$ $1994\}$%
. 
The number of observations for each year above are $\{3874, 3794, 3792,
3666, 4483, 4340,$ $3377\}$. 
In the original data, there are two dependent
variables and 23 variables. 
The two dependent variables are DOCVIS (number
of doctor visits within the last quarter prior to the survey) and HOSPVIS
(number of hospital visits in the last calendar year). 
But in the interest of brevity, we just study on DOCVIS in this paper. 
We list all the variable and its mean and standard deviation in Table S1 in the Supplementary Material. 

Consider the data in each observed year, we build models for DOCVIS by the $\ell_1$ penalized NBR method and classical NBR method via maximum likelihood
estimation. 
We randomly choose 500 samples to train and 500 samples to test for each
year's data. 
For the $\ell_1$ penalized NBR, we use the 10-fold cross validation to
select the penalty level $\lambda$. 
We exhibit the prediction errors (PE) and regression coefficients of both two NBR methods in Table \ref{tab:realdata2}. The "P-NBR" is short for $\ell_1$ penalized NBR. 
The results show that the $\ell_1$ penalized NBR not only simplifies the model but also produces more accurate PE than the classical NBR. 
So we can conjecture that the true model is sparse, in which the most important variables to effect DOCVIS are FEMALE, HSAT and HHKIDS. 
The classical NBR used all variables including some uncorrelated
ones and made a misleading prediction. 

%
%
%

\begin{table}[tbp]
	\centering
	\caption{$\ell_1$-penalized NBR vs general NBR}
	\label{tab:realdata2}%
	\resizebox{6.5in}{!}{
		\begin{tabular}{lccccccccccccccccccccc}
			\hline
			\multirow{2}{*}{Variables} && \multicolumn{2}{c}{1984}&&  \multicolumn{2}{c}{1985}&&  \multicolumn{2}{c}{1986}&&  \multicolumn{2}{c}{1987}&&  \multicolumn{2}{c}{1988}&&  \multicolumn{2}{c}{1991}&&  \multicolumn{2}{c}{1994}\\
			\cline{3-4}\cline{6-7}\cline{9-10}\cline{12-13}\cline{15-16}\cline{18-19}\cline{21-22}
			&& P-NBR & NBR && P-NBR & NBR && P-NBR & NBR && P-NBR & NBR && P-NBR & NBR && P-NBR & NBR && P-NBR & NBR\\
			\hline
			\bf{PE}	&&	44.391	&	2394.968	&&	54.382	&	2905.62	&&	33.194	&	3173.689	&&	75.535	&	1981.128	&&	21.604	&	2156.603	&&	56.350	&	2095.221	&&	39.373	&	2927.628	\\
			Intercept	&&	2.193	&	3.208	&&	1.669	&	0.217	&&	1.945	&	-0.663	&&	2.705	&	3.234	&&	2.053	&	0.333	&&	2.051	&	1.152	&&	2.182	&	-0.053	\\
			FEMALE	&&	0.161	&	0.331	&&	0.066	&	0.265	&&	0.049	&	0.355	&&	0.061	&	0.165	&&	0.039	&	0.246	&&	0.033	&	0.175	&&	0.354	&	0.636	\\
			AGE	&&		&	-0.332	&&	0.395	&	0.731	&&	0.005	&	0.274	&&		&	-0.117	&&		&	-0.244	&&	0.175	&	0.494	&&		&	-0.097	\\
			HSAT	&&	-1.599	&	-1.8	&&	-1.241	&	-1.434	&&	-1.207	&	-1.464	&&	-1.086	&	-1.362	&&	-1.115	&	-1.310	&&	-1.549	&	-1.762	&&	-1.717	&	-2.100	\\
			HANDDUM	&&		&	0.079	&&		&	-0.021	&&	0.182	&	0.168	&&	-0.214	&	-0.455	&&	0.032	&	0.052	&&	0.010	&	0.092	&&	0.092	&	0.577	\\
			HANDPER	&&		&	0.048	&&	0.091	&	0.148	&&		&	0.054	&&		&	-0.041	&&	0.092	&	0.142	&&		&	-0.060	&&		&	-0.479	\\
			HHNINC	&&		&	0.019	&&		&	-0.177	&&		&	-0.162	&&	-0.088	&	-0.206	&&	-0.181	&	-0.489	&&		&	0.207	&&	-0.050	&	-0.462	\\
			HHKIDS	&&	-0.028	&	-0.214	&&	-0.024	&	-0.128	&&	-0.031	&	-0.159	&&	-0.038	&	-0.136	&&		&	-0.167	&&	-0.007	&	-0.088	&&		&	-0.040	\\
			EDUC	&&		&	-0.344	&&		&	1.425	&&		&	2.872	&&	-0.111	&	0.028	&&		&	2.123	&&		&	0.407	&&		&	1.700	\\
			MARRIED	&&		&	0.014	&&		&	0.24	&&		&	0.000	&&	-0.081	&	-0.210	&&		&	0.140	&&		&	-0.057	&&		&	-0.059	\\
			HAUPTS	&&		&	0.13	&&		&	-0.257	&&		&	-0.394	&&		&	-0.054	&&	0.006	&	0.431	&&		&	-0.075	&&		&	0.960	\\
			REALS	&&		&	0.082	&&		&	-0.102	&&		&	-0.334	&&		&	0.043	&&		&	0.128	&&		&	0.070	&&		&	0.445	\\
			FACHHS	&&		&	0.058	&&		&	-0.165	&&		&	-0.270	&&	-0.014	&	-0.156	&&		&	-0.061	&&		&	-0.002	&&		&	0.092	\\
			ABITUR	&&		&	-0.039	&&		&	-0.176	&&		&	-0.549	&&		&	0.022	&&	-0.040	&	-0.200	&&		&	0.016	&&		&	0.136	\\
			UNIV	&&		&	0.152	&&		&	-0.217	&&	-0.002	&	-0.319	&&		&	-0.160	&&		&	-0.179	&&		&	-0.040	&&		&	-0.026	\\
			WORKING	&&		&	0.155	&&		&	0.327	&&		&	-0.256	&&	-0.214	&	-0.038	&&		&	-0.448	&&	-0.025	&	0.112	&&		&	-0.079	\\
			BLUES	&&		&	-0.104	&&		&	-0.117	&&		&	0.309	&&		&	-0.131	&&		&	0.369	&&		&	-0.058	&&		&	0.049	\\
			WHITEC	&&		&	-0.092	&&		&	-0.155	&&		&	0.291	&&		&	-0.111	&&	-0.006	&	0.277	&&	-0.035	&	-0.313	&&		&	0.125	\\
			SELF	&&		&	-0.129	&&	-0.02	&	-0.219	&&	0.016	&	0.238	&&	-0.063	&	-0.212	&&		&	0.135	&&		&	-0.087	&&		&	0.005	\\
			BEAMT	&&		&	-0.228	&&		&	-0.171	&&		&	0.267	&&		&	-0.006	&&		&	0.184	&&		&	-0.051	&&		&	0.107	\\
			PUBLIC	&&		&	-0.407	&&		&	0.095	&&	0.177	&	0.520	&&	0.080	&	0.328	&&		&	-0.240	&&		&	0.297	&&		&	0.107	\\
			ADDON	&&		&	-0.149	&&		&	0.02	&&		&	-0.110	&&		&	0.008	&&		&	-0.020	&&		&	-0.067	&&		&	-0.044	\\
			\hline
	\end{tabular}}
\end{table}

\section{Conclusion}\label{sec:conclusion}
We studied on two theoretical choices of penalty level, under which we proved that the $\ell_1$ penalized negative binomial estimator is $\ell_1$ consistent. 
We then conduct a simulation, whose results further confirm the convergence tendency of estimation errors. 
Finally, an real application shows that the $\ell_1$ penalized NBR can produce a concise model that only contains key variables, and it has much smaller prediction errors than the classical NBR.


\appendix \section{Proof of Theorem 2.1}
Let $\difvec = \est - \target$. Recall that $\nzeroindx = \{j:\targets_{j} \neq 0\}$. By definition of 
$\est$ and the convexity of $L(\parameter)$, we have 
\begin{equation}  \label{ineq:bydefinition}
\begin{split}
L(\est)-L(\target)&\le\lambda(\|\target\|_{1}-\|\est\|_{1})
\\
&= \lambda [(\|\target_{\nzeroindx} \|_1 - \|\est_{\nzeroindx} \|_1) +(\|\target_{\nzeroindx^{c}}\| -
\|\est_{\nzeroindx^{c}}\|_1)] \\
&\le\lambda(\|\difvec_{\nzeroindx}\|_{1}-\|\difvec_{\nzeroindx^{c}}\|_{1}),
\end{split}%
\end{equation}
and 
\begin{equation}  \label{ineq:convexity}
L(\est)-L(\target)\ge\difvec\tp \nabla
L(\target)\ge-V\|\difvec\|_{1}\ge-\frac{\lambda}{c}\|\difvec\|_{1},
\end{equation}
where $V=\Vert \nabla L(\target)\Vert _{\infty }
$, and the last inequality utilizes the condition $\lambda > cV$. Combining (%
\ref{ineq:bydefinition}) and (\ref{ineq:convexity}), we obtain that 
\begin{equation*}
\|\difvec_{\nzeroindx^{c}}\|_{1}\le \frac{c+1}{c-1}\|\difvec_{\nzeroindx}\|_{1},
\end{equation*}
which makes the condition \textbf{C4} hold with $\gamma=\frac{c+1}{c-1}>1$.%

\indent For any $\parameter,\vecu,\vecv\in\mathbb{R}^p$, we have 
\begin{equation*}
\begin{split}
\nabla^2 L(\parameter)[\vecv,\vecv]&=\frac{1}{n}\sum%
\limits_{i=1}^{n}(\designrow_{i}\tp\vecv)^{2}\frac{r\exp{\{\designrow_{i}\tp\parameter\}}(\outcomes_{i}+r)}{(r+\exp{\{\designrow_{i}\tp\parameter\}})^2} \\
\nabla^3 L(\parameter)[\vecu,\vecv,\vecv]&=\frac{1}{n}\sum%
\limits_{i=1}^{n}\designrow_{i}\tp\vecu(\designrow_{i}\tp\vecv)^{2}\frac{r\exp{\{\designrow_{i}\tp\parameter\}}(\outcomes_{i}+r)(r-\exp{\{\designrow_{i}\tp\parameter\}})}{(r+\exp{\{\designrow_{i}\tp\parameter\}})^3}.
\end{split}%
\end{equation*}

Under the assumption (C2), it is easy to verify that 
\begin{equation*}
\begin{split}
|\nabla^3 L(\parameter)[\vecu,\vecv,\vecv]|&\le \sup\limits_{i\in[n]}|\designrow_{i}\tp\vecu|\nabla^2
L(\parameter)[\vecv,\vecv] \\
&\le \sup\limits_{i\in[n],j\in[p]}|x_{ij}|\|\vecu\|_{1}\nabla^2 L(\parameter)[\vecv,\vecv] \\
&\le \contx\|\vecu\|_{1}\nabla^2 L(\parameter)[\vecv,\vecv].
\end{split}%
\end{equation*}

Setting $\vecu=\difvec=\est-\target$, we have 
\begin{equation}  \label{ineq:self-corcondient}
\begin{split}
|\nabla^3 L(\parameter)[\vecu,\vecv,\vecv]|&\le \contx(1+ \frac{c+1}{c-1})\|\difvec_{\nzeroindx}\|_{1}%
\nabla^2 L(\parameter)[\vecv,\vecv] \\
&\le \frac{2c\contx\sqrt{s}}{c-1}\|\difvec_{\nzeroindx}\|_{2}\nabla^2 L(\parameter)[\vecv,\vecv].
\end{split}%
\end{equation}
Denoting $\tilde{\contx}=\frac{2c\contx\sqrt{s}}{c-1}$, (\ref{ineq:self-corcondient})
becomes $|\nabla^3 L(\parameter)[\vecu,\vecv,\vecv]|\le \tilde{\contx}\|\difvec_{\nzeroindx}\|_{2}\nabla^2
L(\parameter)[\vecv,\vecv].$

Thus, $L(\cdot)$ is a self-concordant like function with parameter $\tilde{\contx}
$ with respect to $l_2$-norm.By Theorem 6.1 of [13] and condition 
\textbf{C4}, we have 
\begin{equation}  \label{ineq:taylor}
\begin{split}
L(\est)-L(\target)&\ge\difvec\tp \nabla L(\target)+\frac{%
	\difvec\tp \nabla^2 L(\target)\difvec}{\tilde{\contx}^{2}\|\difvec_{\nzeroindx}\|_{2}^{2}}%
(\exp{\{-\tilde{\contx}\|\difvec_{\nzeroindx}\|_{2}\}}+\tilde{\contx}\|\difvec_{\nzeroindx}\|_{2}-1) \\
&\ge-\|\nabla L(\target)\|_{\infty}\|\difvec\|_{1}+\frac{\difvec\tp \nabla^2
	L(\target)\difvec}{\tilde{\contx}^{2}\|\difvec_{\nzeroindx}\|_{2}^{2}}(\exp{\{-\tilde{\contx}%
	\|\difvec_{\nzeroindx}\|_{2}\}}+\tilde{\contx}\|\difvec_{\nzeroindx}\|_{2}-1) \\
&\ge-\frac{\lambda}{c}\|\difvec\|_{1}+\frac{\difvec\tp \nabla^2
	L(\target)\difvec}{\tilde{\contx}^{2}\|\difvec_{\nzeroindx}\|_{2}^{2}}(\exp{\{-\tilde{\contx}%
	\|\difvec_{\nzeroindx}\|_{2}\}}+\tilde{\contx}\|\difvec_{\nzeroindx}\|_{2}-1)
\end{split}%
\end{equation}
Combining (\ref{ineq:bydefinition}) and (\ref{ineq:taylor}), we have 
\begin{equation}  \label{ineq:1}
\begin{split}
&\quad\frac{\difvec\tp \nabla^2 L(\target)\difvec}{\tilde{\contx}%
	^{2}\|\difvec_{\nzeroindx}\|_{2}^{2}}(\exp{\{-\tilde{\contx}\|\difvec_{\nzeroindx}\|_{2}\}}+\tilde{\contx}%
\|\difvec_{\nzeroindx}\|_{2}-1) \\
&\le\lambda\|\difvec_{\nzeroindx} \|_1+\frac{\lambda}{c}\|\difvec\|_{1}\le \frac{c+1}{c-1}%
\lambda\|\difvec_{\nzeroindx}\|_{1}\le \frac{c+1}{c-1}\lambda\sqrt{s}%
\|\difvec_{\nzeroindx}\|_{2}.
\end{split}%
\end{equation}
By condition \textbf{C4} and (\ref{ineq:1}), we have 
\begin{equation}  \label{ineq:2}
\exp{\{-\tilde{\contx}\|\difvec_{\nzeroindx}\|_{2}\}}+\tilde{\contx}\|\difvec_{\nzeroindx}\|_{2}-1\le\frac{%
	(c+1)\lambda \sqrt{s}\tilde{\contx}^{2}}{(c-1)\phi_0^{2}}\|\difvec_{\nzeroindx}\|_{2}.
\end{equation}
Set 
\begin{equation}  \label{eqn:h}
h=\frac{(c+1)\lambda \sqrt{s}\tilde{\contx}}{(c-1)\phi_0^{2}} = \frac{%
	2c(c+1)\lambda \contx s}{(c-1)^{2}\phi_0^2},
\end{equation}
then according to the condition on $\lambda$ such that $\lambda s\le\frac{%
	(c-1)^{2}\phi_0^{2}}{6c\contx (c+1)}$, we have $h\le\frac1{3}$. Denote $w=\tilde{\contx}%
\|\difvec_{\nzeroindx}\|_{2}$, then to solve (\ref{ineq:2}) is equivalent to solve the
inequality $\exp{\{-w\}}+w-1\le hw$. By Taylor formula, we have $\frac{w^2}{2}-%
\frac{w^3}{6}\le \exp{\{-w\}}+w-1\le hw$ which implies $\{w:\exp{\{-w\}}+w-1\le hw, h\le
\frac1{3}\}\subseteq\{w:\frac{w^2}{2}-\frac{w^3}{6}\le hw, h\le \frac1{3}\}$%
. Since under the condition $h\le \frac1{3}$, the solution of inequality $%
\frac{w^2}{2}-\frac{w^3}{6}\le hw$ is $w\le c_{1}h$ for some constant $%
c_{1}\in(2,3]$, then 
\begin{equation*}
\{w:\exp{\{-w\}}+w-1\le hw, h\le \frac1{3}\}\subseteq\{w\le c_{1}h \mathrm{\ with\
	some\ constant\ }c_{1}\in(2,3]\}.
\end{equation*}
So, from (\ref{ineq:2}), we obtain 
\begin{equation*}
\tilde{\contx}\|\difvec_{\nzeroindx}\|_{2}\le\frac{c_{1}\lambda \sqrt{s}(c+1)}{%
	\phi_0^{2}(c-1)}\tilde{\contx},
\end{equation*}
that is, 
\begin{equation}  \label{ineq:4}
\|\difvec_{\nzeroindx}\|_{2}\le\frac{c_{1}\lambda \sqrt{s}(c+1)}{\phi_0^{2}(c-1)}.
\end{equation}
Hence, notice the relationship $\|\difvec\|_{1}\le(1+\frac{c+1}{c-1})\sqrt{s}%
\|\difvec_{\nzeroindx}\|_{2}=\frac{2c\sqrt{s}}{c-1}\|\difvec_{\nzeroindx}\|_{2}$, by (\ref%
{ineq:4}) we have 
\begin{equation}  \label{ineq:5}
\|\difvec\|_{1}\le\frac{2c_{1}c(c+1)}{(c-1)^{2}\phi_0^{2}}\lambda s.
\end{equation}
Then, (2.2) is obtained. Furthermore, by (\ref%
{ineq:bydefinition}) and (\ref{ineq:convexity}), we obtain 
\begin{equation}  \label{ineq:6}
|L(\est)-L(\target)|\le\lambda\|\difvec\|_1\le\frac{2c_{1}c(c+1)}{%
	(c-1)^{2}\phi_0^{2}}\lambda^{2} s,
\end{equation}
which implies (2.3). We finish the proof. $\square$

\section{Proof of Lemma 2.3}
Recall the gradient of $L(\parameter)$ at the point $\parameter=\target$ 
\begin{equation*}
\nabla L(\target) = -\frac1{n}\sum\limits_{i=1}^{n}\designrow_{i}\times v_i\times 
\frac{(\outcomes_{i}-\exp{\{\designrow_{i}\tp \target\}})}{\sqrt{\frac{\exp{\{\designrow_{i}\tp %
				\target\}}(r+\exp{\{\designrow_{i}\tp \target\}})}{r}}},
\end{equation*}
where 
\begin{equation*}
v_{i}=\sqrt{\frac{r\exp{\{\designrow_{i}\tp \target\}}}{\exp{\{\designrow_{i}\tp \target\}}+r}}. 
\end{equation*}%
Denote $\epsilon_{i}=(\outcomes_{i}-\exp{\{\designrow_{i}\tp \target\}})/\sqrt{\exp{\{\designrow_{i}\tp \target\}}(r+\exp{\{\designrow_{i}\tp \target\}})/r}$, then 
\begin{equation*}
V= \|\nabla L(\target) \|_{\infty}=\max\limits_{j\in[p]}|\frac1{n}\sum%
\limits_{i=1}^{n}x_{ij}v_i\epsilon_{i}|.
\end{equation*}
Recall $v_{n}=\max\limits_{1\leq i\leq
	n}v_{i}$ and denote $t_{p,\alpha}=\Phi^{-1}(1-\frac{\alpha}{2p})$, then $\lambda= cv_n(%
\sqrt{n})^{-1}t_{p,\alpha}$. Hence 
\begin{equation}  \label{ineq:prob}
\begin{split}
\mathbb{P}(cV>\lambda)&=\mathbb{P}(\max\limits_{j\in[p]}|\frac1{n}\sum%
\limits_{i=1}^{n}x_{ij}v_i\epsilon_{i}|>v_n(\sqrt{n})^{-1}t_{p,\alpha}) \\
&\le \mathbb{P}(\max\limits_{j\in[p]}|\frac1{n}\sum\limits_{i=1}^{n}x_{ij}%
\epsilon_{i}|>(\sqrt{n})^{-1}t_{p,\alpha}) \\
&\le p\max\limits_{j\in[p]}\mathbb{P}(|\sum\limits_{i=1}^{n}x_{ij}%
\epsilon_{i}|>\sqrt{n}t_{p,\alpha}).
\end{split}%
\end{equation}
Since $\outcomes_{i}|\designrow_{i}\sim \mathrm{NB}(r,\mu_i)$ with $\mu_{i}=\exp{\{\designrow_{i}\tp %
	\target\}}$, then 
\begin{equation*}
\mathbb{E}(\exp{\left\{\theta\epsilon_{i}\right\}})=\exp{\left\{-\theta\sqrt{\frac{r\mu_i}{\mu_i+r}}%
	\right\}}\left(1+\frac{\mu_i}{r}\left(1-\exp{\left\{\theta\sqrt{\frac{r}{\mu_i(\mu_i+r)}}%
	\right\}}\right)\right)^{-r}
\end{equation*}
is a positive constant for all $\theta<\sqrt{\frac{\mu_i(\mu_i+r)}{r}}\ln(%
\frac{r}{\mu_i}+1)$. By the exponential Chebyshev's inequality, we have 
\begin{equation}  \label{ineq:sub-exponential}
\mathbb{P}(|\epsilon_{i}|>A)<\exp{\{-A/w_{1}\}}\mathbb{E}(\exp{%
	\{\epsilon_{i}/w_{1}\}})=C_{1}\exp{\{-A/w_{1}\}}
\end{equation}
with some constant $C_{1}=\mathbb{E}(\exp{\{\epsilon_{i}/w_{1}\}})>0$ and $%
w_{1}>\left(\sqrt{\frac{\mu_i(\mu_i+r)}{r}}\ln(\frac{r}{\mu_i}%
+1)\right)^{-1} $. Denote $\hat{\epsilon}_{i}=\epsilon_{i}1_{\{|%
	\epsilon_{i}|\le A\}}$ and $\check{\epsilon}_{i}=\epsilon_{i}1_{\{|%
	\epsilon_{i}| > A\}}$. Taking $A=3w_{1}\log{p}$, we have 
\begin{equation*}
\begin{split}
\mathbb{P}(|\sum\limits_{i=1}^{n}x_{ij}\epsilon_{i}|>\sqrt{n}t_{p,\alpha})&= 
\mathbb{P}(|\sum\limits_{i=1}^{n}x_{ij}(\hat{\epsilon}_{i}+\check{\epsilon}%
_{i})|>\sqrt{n}t_{p,\alpha},\sup\limits_{i\in[n]}|\epsilon_{i}|\le A) \\
&\quad+\mathbb{P}(|\sum\limits_{i=1}^{n}x_{ij}(\hat{\epsilon}_{i}+\check{%
	\epsilon}_{i})|>\sqrt{n}t_{p,\alpha},\sup\limits_{i\in[n]}|\epsilon_{i}|> A)
\\
&\le\mathbb{P}(|\sum\limits_{i=1}^{n}x_{ij}\hat{\epsilon}_{i}|>\sqrt{n}%
t_{p,\alpha})+\mathbb{P}(\sup\limits_{i\in[n]}|\epsilon_{i}|> A).
\end{split}%
\end{equation*}

Denote $P_{1}=\mathbb{P}(|\sum\limits_{i=1}^{n}x_{ij}\hat{\epsilon}_{i}|>%
\sqrt{n}t_{p,\alpha})$ and $P_{2}=\mathbb{P}(\sup\limits_{i\in[n]%
}|\epsilon_{i}|> A)$, then the above inequality can be written as 
\begin{equation}  \label{ineq:truncation}
\mathbb{P}(|\sum\limits_{i=1}^{n}x_{ij}\epsilon_{i}|>\sqrt{n}%
t_{p,\alpha})\le P_{1}+P_{2}.
\end{equation}

By inequality (\ref{ineq:sub-exponential}) with $A=3w_{1}\log{p}$, we obtain
that 
\begin{equation}  \label{ineq:P2}
P_{2}\le\sum\limits_{i=1}^{n}\mathbb{P}(|\epsilon_{i}|> A)\le C_{1}n\exp{\{-3\log%
	{p}\}}=C_{1}n/p^{3}.
\end{equation}

To estimate the $P_{1}$, we need the following Sakhanenko type moderate
deviation theorem of \citep{Sak91}, i.e.

\begin{lem}
	\label{lem:MDT} Let $\eta_{1},\cdots,\eta_{n}$ be independent random
	variables with $\mathbb{E}\eta_{i}=0$ and $|\eta_{i}|<1$ for all $i\in[n]$.
	Denote $\sigma_{n}^{2}=\sum\limits_{i=1}^{n}\mathbb{E}\eta_{i}^{2}$ and $%
	T_{n}=\sum\limits_{i=1}^{n}\mathbb{E}|\eta_{i}|^{3}/\sigma_{n}^{3}$. Then
	there exists a positive constant $D$ such that for all $x\in[%
	1,\frac1{D}\min\{\sigma_{n},L_{n}^{-1/3}\}]$ 
	\begin{equation*}
	\mathbb{P}(\sum\limits_{i=1}^{n}\eta_{i}>x\sigma_{n})=(1+O(1)x^{3}T_{n})\bar%
	\Phi(x),
	\end{equation*}
	where $\bar\Phi(x)=1-\Phi(x)$ and $\Phi(x)$ is the cumulative distribution
	function of standard normal distribution.
\end{lem}

Since $\mathbb{E}(\epsilon_{i})=\mathbb{E}(\hat{\epsilon}_{i})+\mathbb{E}(%
\check{\epsilon}_{i})=0$, then it is easy to obtain that 
\begin{align*}
|\mathbb{E} \hat{\epsilon}_{i}|&=|\mathbb{E} \check{\epsilon}_{i}|\le\mathbb{%
	E}|\check{\epsilon}_{i}|=\mathbb{E}|\epsilon_{i}|1_{\{|\epsilon_{i}|>A\}}=%
\int_{A}^{+\infty}z d F(z)+\int_{-\infty}^{-A}-z d F(z) \\
&=\left\{z(F(z)-1)|_{A}^{+\infty}-\int_{A}^{+\infty}(F(z)-1)dz\right\}+\left%
\{\int_{-\infty}^{-A}F(z)dz-zF(z)|_{-\infty}^{-A}\right\} \\
&\le A(1-F(A))+\int_{A}^{+\infty}C_{1}\exp{\{-z/w_{1}\}}dz +
\int_{-\infty}^{-A}C_{1}\exp{\{z/w_{1}\}}dz+AF(-A) \\
&\le C_{1}(A+2w_{1})\exp{\{-A/w_{1}\}}\le 2C_{1}A\exp{\{-A/w_{1}\}}=\frac{6C_{1}w_{1}\log%
	{p}}{p^{3}},
\end{align*}
where the last second and third inequalities utilize the relations $F(a)=%
\mathbb{P}(\epsilon_{i}>a)\le C_{1}\exp{\{-a/w_{1}\}}$ and $F(-a)=\mathbb{P}%
(\epsilon_{i}<-a)\le C_{1}\exp{\{-a/w_{1}\}}$ for any $a>0$. Denote $m=6C_{1}w_{1}\log{p}/p^{3}$, then $|\mathbb{E} \hat{\epsilon}_{i}|\le m$ and 
$m=o(n^{-2})$. 

Since 
\begin{equation}  \label{ineq:P1-1}
\begin{split}
P_{1}&=\mathbb{P}(|\sum\limits_{i=1}^{n}x_{ij}(\hat{\epsilon}_{i}-\mathbb{E}%
\hat{\epsilon}_{i}+\mathbb{E}\hat{\epsilon}_{i})|>\sqrt{n}t_{p,\alpha}) \\
&\le\mathbb{P}(|\sum\limits_{i=1}^{n}x_{ij}(\hat{\epsilon}_{i}-\mathbb{E}%
\hat{\epsilon}_{i}))|>\sqrt{n}t_{p,\alpha}-|\sum\limits_{i=1}^{n}x_{ij}%
\mathbb{E}\hat{\epsilon}_{i}|),
\end{split}%
\end{equation}
we need to estimate $|x_{ij}(\hat{\epsilon}_{i}-\mathbb{E} \hat{\epsilon}%
_{i})$ and $|\sum\limits_{i=1}^{n}x_{ij}\mathbb{E}\hat{\epsilon}_{i}|$. By
condition \textbf{C2}, 
\begin{equation*}
|x_{ij}(\hat{\epsilon}_{i}-\mathbb{E}\hat{\epsilon}_{i})|\le(\sup\limits_{i%
	\in[n],j\in[p]}|x_{ij}|)(|\hat{\epsilon}_{i}|+|\mathbb{E}\hat{\epsilon}%
_{i}|)\le \contx(A+m).
\end{equation*}
By Cauchy-Schwarz inequality, 
\begin{equation*}
|\sum\limits_{i=1}^{n}x_{ij}\mathbb{E}\hat{\epsilon}_{i}|\le\sqrt{%
	(\sum\limits_{i=1}^{n}x_{ij}^{2})(\sum\limits_{i=1}^{n}|\mathbb{E}\hat{%
		\epsilon}_{i}|^{2})}\le nm.
\end{equation*}

Denoting $\eta_{ij}=x_{ij}(\hat{\epsilon}_{i}-\mathbb{E}\hat{\epsilon}%
_{i})/\contx(A+m)$, we have $\mathbb{E} \eta_{ij}=0$ and $|\eta_{ij}|<1$. Notice
that $\mathbb{E} \hat{\epsilon}_{i}^{2}\le \mathbb{E} \epsilon_{i}^{2}=1$.
Denoting $\sigma_{nj}^{2}=\sum\limits_{i=1}^{n}\mathbb{E}\eta_{ij}^{2}$ and $%
T_{nj}=\sum\limits_{i=1}^{n}\mathbb{E}|\eta_{ij}|^{3}/\sigma_{nj}^{3}$, we
have

\begin{align*}
\sigma_{nj}^{2}&=\frac1{\contx^2(A+m)^{2}}\sum\limits_{i=1}^{n}\mathbb{E}%
(x_{ij}^{2}(\hat{\epsilon}_{i}-\mathbb{E}\hat{\epsilon}_{i})^{2}) \\
&\le\frac1{\contx^2(A+m)^{2}}\sum\limits_{i=1}^{n}x^{2}_{ij}\mathbb{E}\hat{%
	\epsilon}_{i}^{2}\le\frac1{\contx^2(A+m)^{2}}\sum\limits_{i=1}^{n}x_{ij}^{2} \\
&=\frac{n}{\contx^2(A+m)^{2}}, \\
T_{nj}&\le \sum\limits_{i=1}^{n}\mathbb{E}|\eta_{ij}|^{2}/\sigma_{nj}^{3}=%
\frac1{\sigma_{nj}}.
\end{align*}
Hence, $\sigma_{nj}^{2}=O(\frac{n}{(A+m)^2})$ and $L_{nj}=O(\frac{A+m}{\sqrt{%
		n}})$. By inequality \eqref{ineq:P1-1} and Lemma \ref{lem:MDT}, for large
enough $n,p$ such that $n\le p \le o(\exp{\{n^{1/5}\}})$(condition \textbf{C3}),
we have 
\begin{equation}  \label{ineq:P1-2}
\begin{split}
P_{1}&\le \mathbb{P}(|\sum\limits_{i=1}^{n}\frac{x_{ij}(\hat{\epsilon}_{i}-%
	\mathbb{E}\hat{\epsilon}_{i}))}{\contx(A+m)}|>\frac{\sqrt{n}}{\contx(A+m)}%
(t_{p,\alpha}-\sqrt{n}m)) \\
&\le \mathbb{P}(|\sum\limits_{i=1}^{n}\eta_{ij}|>\sigma_{nj}(t_{p,\alpha}-%
\sqrt{n}m)) \\
&=2\left(1+O(1)\left(t_{p,\alpha}-\sqrt{n}m\right)^{3}T_{nj}\right)\bar{\Phi}%
(t_{p,\alpha}-\sqrt{n}m) \\
\end{split}%
\end{equation}
with $t_{p,\alpha}-\sqrt{n}m$ uniformly in $[1, O(n^{1/6}(\log{p})^{-1/3})]$%
. Next, we estimate $O(1)(t_{p,\alpha}-\sqrt{n}m)^{3}T_{nj}$ and $\bar{\Phi}%
(t_{p,\alpha}-\sqrt{n}m)$ respectively. Notice that $\log{(p/\alpha)}%
<t_{p,\alpha}^{2}<2\log{(2p/\alpha)}$ when $p/\alpha>8$. Then, under
condition \textbf{C3}, we have 
\begin{equation}  \label{ineq:o(1)}
\begin{split}
&\quad O(1)(t_{p,\alpha}-\sqrt{n}m)^{3}T_{nj} \\
&= O(1)(\sqrt{2\log{(2p/\alpha)}}-\sqrt{n}m)^{3}n^{-1/2}(3w_{1}\log{p}+m).
\end{split}%
\end{equation}
Furthermore, by the fact that for all $a>0$ the inequality $\frac{a}{1+a^{2}}%
\phi(a)\le\bar{\Phi}(a)\le\frac{\phi(a)}{a}$ holds where $\phi(\cdot)$ is
the density function of standard normal distribution, we have 
\begin{equation}  \label{ineq:Phibar}
\begin{split}
\bar{\Phi}(t_{p,\alpha}-\sqrt{n}m)&\le\frac{\phi(t_{p,\alpha}-\sqrt{n}m)}{%
	t_{p,\alpha}-\sqrt{n}m}=\phi(t_{p,\alpha})\frac{\exp\{t_{p,\alpha}\sqrt{n}%
	m-nm^{2}/2\}}{t_{p,\alpha}-\sqrt{n}m} \\
&=\frac{t_{p,\alpha}}{1+t_{p,\alpha}^{2}}\phi(t_{p,\alpha})\frac{%
	1+t_{p,\alpha}^{2}}{t_{p,\alpha}(t_{p,\alpha}-\sqrt{n}m)}\exp\{t_{p,\alpha}%
\sqrt{n}m-nm^{2}/2\} \\
&\le\bar{\Phi}(t_{p,\alpha})\frac{1+t_{p,\alpha}^{2}}{t_{p,\alpha}(t_{p,%
		\alpha}-\sqrt{n}m)}\exp\{t_{p,\alpha}\sqrt{n}m-nm^{2}/2\} \\
&= \frac{\alpha}{2p}(1+\frac1{t_{p,\alpha}^{2}})\frac1{1-\sqrt{n}%
	m/t_{p,\alpha}}\exp\{t_{p,\alpha}\sqrt{n}m-nm^{2}/2\} \\
&\le\frac{\alpha}{2p}(1+\frac1{\log{(p/\alpha)}})\frac{\exp\{(n\log{%
		(p/\alpha)})^{1/2}m-nm^{2}/2\}}{1-\sqrt{n}m/(\log{(p/\alpha)})^{1/2}}.
\end{split}%
\end{equation}
Combining (\ref{ineq:P1-2}), (\ref{ineq:o(1)}) and (\ref{ineq:Phibar}), we
have 
\begin{equation}  \label{ineq:P1}
\begin{split}
P_{1}&\le\frac{\alpha}{p}(1+O(1)(\sqrt{2\log{(2p/\alpha)}}-\sqrt{n}%
m)^{3}n^{-1/2}(3w_{1}\log{p}+m)) \\
&\quad\times(1+\frac1{\log{(p/\alpha)}})\frac{\exp\{(n\log{(p/\alpha)}%
	)^{1/2}m-nm^{2}/2\}}{1-\sqrt{n}m/(\log{(p/\alpha)})^{1/2}}.
\end{split}%
\end{equation}
Hence, combining (\ref{ineq:prob}), \eqref{ineq:truncation}, (\ref{ineq:P2})
and (\ref{ineq:P1}), we have 
\begin{equation*}
\begin{split}
\mathbb{P}(\lambda<cV)&\le p(P_{1}+P_{2}) \\
&\le\alpha(1+O(1)(\sqrt{2\log{(2p/\alpha)}}-\sqrt{n}m)^{3}n^{-1/2}(3w_{1}\log%
{p}+m)) \\
&\quad\times(1+\frac1{\log{(p/\alpha)}})\frac{\exp\{(n\log{(p/\alpha)}%
	)^{1/2}m-nm^{2}/2\}}{1-\sqrt{n}m/(\log{(p/\alpha)})^{1/2}}+C_{1}n/p^{2},
\end{split}%
\end{equation*}
where $C_{1}$ and $w_{1}$ are some positive constants. 

\indent So, the probability of event $\{\lambda\ge cV\}$ is 
\begin{equation*}
\begin{split}
\mathbb{P}(\lambda\ge cV)&\ge 1 - \alpha(1+O(1)(\sqrt{2\log{(2p/\alpha)}}-%
\sqrt{n}m)^{3}n^{-1/2}(3w_{1}\log{p}+m)) \\
&\qquad\quad\times(1+\frac1{\log{(p/\alpha)}})\frac{\exp\{(n\log{(p/\alpha)}%
	)^{1/2}m-nm^{2}/2\}}{1-\sqrt{n}m/(\log{(p/\alpha)})^{1/2}}-C_{1}n/p^{2}.
\end{split}%
\end{equation*}
Additionally, notice that $m$, $\sqrt{n}m$ and $nm^{2}$ are $o(n^{-2})$. As $%
n,p\rightarrow\infty$ with $n\le p \le o(\exp{\{n^{1/5}\}})$, it is easy to obtain
that 
\begin{equation*}
\mathbb{P}(\lambda\ge cV)\le 1-\alpha(1+o(1)).
\end{equation*}


\section{Table}
\begin{table}[h]
	\centering
	\caption{list of variables}
	\label{tab:realdata1}%
	\resizebox{5.8in}{!}{
		\begin{tabular}{lclcccc}
			\hline
			Variables && Description && Mean && SD \\
			\hline
			DOCVIS    && number of doctor visits in last three months      &&  3.184	&&	5.69
			\\
			HOSPVIS   && number of hospital visits in last calendar year   &&  0.138	&&	0.884
			\\
			ID        && person - identification number, $1,\cdots,7293$     &&     	&&	
			\\
			FEMALE    && female = 1; male = 0               &&  0.479	&&	0.5
			\\
			1984      && Year = 1984 (0/1)     &&  0.142	&&	0.349
			\\
			1985      && Year = 1985 (0/1)      &&  0.139	&&	0.346
			\\
			1986      && Year = 1986 (0/1)      &&  0.139	&&	0.346
			\\
			1987      && Year = 1987 (0/1)      &&  0.134	&&	0.341
			\\
			1988      && Year = 1988 (0/1)      &&  0.164	&&	0.37
			\\
			1991      && Year = 1991 (0/1)      &&  0.159	&&	0.366
			\\
			1994      && Year = 1994 (0/1)      &&  0.124	&&	0.329
			\\
			
			AGE       && age in years    &&  43.526	&&	11.33
			\\
			HSAT      && health satisfaction, coded 0 (low) - 10 (high)    &&  6.785	&&	2.294
			\\
			HANDDUM   && handicapped = 1; otherwise = 0    &&  0.214	&&	0.41
			\\
			HANDPER   && degree of handicap in percent (0 - 100)     &&  7.012	&&	19.265
			\\
			HHNINC    && household nominal monthly net income in German marks / 1000     &&  3.521	&&	1.769
			\\
			HHKIDS    && children under age 16 in the household = 1; otherwise = 0     &&  0.403	&&	0.49
			\\
			EDUC      && years of schooling     &&  11.321	&&	2.325
			\\
			MARRIED   && married = 1; otherwise = 0      &&  0.759	&&	0.428
			\\
			HAUPTS    && highest schooling degree is Hauptschul degree = 1; otherwise =
			0     &&  0.624	&&	0.484
			\\
			REALS     && highest schooling degree is Realschul degree = 1; otherwise = 0      &&  0.197	&&	0.398
			\\
			FACHHS    && highest schooling degree is Polytechnical degree = 1;
			otherwise = 0        &&  0.041	&&	0.198
			\\
			ABITUR    && highest schooling degree is Abitur = 1; otherwise = 0         &&  0.117	&&	0.321
			\\
			UNIV      && highest schooling degree is university degree = 1; otherwise = 0      &&  0.072	&&	0.258
			\\
			WORKING   && employed = 1; otherwise = 0      &&  0.677	&&	0.468
			\\
			BLUEC     && blue collar employee = 1; otherwise = 0       &&  0.244	&&	0.429
			\\
			WHITEC    && white collar employee = 1; otherwise = 0      &&  0.3	&&	0.458
			\\
			SELF      && self employed = 1; otherwise = 0      &&  0.062	&&	0.241
			\\
			BEAMT     && civil servant = 1; otherwise = 0      &&  0.075	&&	0.263
			\\
			PUBLIC    && insured in public health insurance = 1; otherwise = 0      &&  0.886	&&	0.318
			\\
			ADDON     && insured by add-on insurance = 1; otherswise = 0      &&  0.019	&&	0.136
			\\
			\hline
	\end{tabular}}
\end{table}





\bibliographystyle{elsarticle-harv} 
\bibliography{elsarticle-template-num-names}





\end{document}